\documentclass[12pt]{article}

\usepackage{amsmath,amsthm,amsfonts}[1996/11/01]

\advance \textheight by -1 \topmargin
\advance \textheight by -1 \headheight
\advance \textheight by -1 \headsep
\advance \textheight by -2 \footskip
\vsize = \textheight
\hsize = \textwidth
\DeclareMathOperator{\Aut}{Aut}

\DeclareMathOperator{\GL}{GL}
\DeclareMathOperator{\BW}{BW}

\theoremstyle{plain}
\newtheorem{theorem}{Theorem}
\newtheorem{lemma}[theorem]{Lemma}

\newtheorem{proposition}[theorem]{Proposition}

\newtheorem{definition}[theorem]{Definition}
\theoremstyle{remark}

\newtheorem{remark}[theorem]{Remark}

\numberwithin{theorem}{section}

\newcommand{\bew}{\noindent\underline{Proof.}\ }

\newcommand{\disj}{\stackrel{.}{\cup}}

\newcommand{\Z}{{\mathbb{Z}}}

\newcommand{\F}{{\mathbb{F}}}

\newcommand{\N}{{\mathbb{N}}}
\newcommand{\R}{{\mathbb{R}}}

\newcommand{\eb}{\phantom{zzz}\hfill{$\square $}\smallskip}

\newcommand{\cV}{{\mathcal V}}
\newcommand{\cU}{{\mathcal U}}
\renewcommand{\em}{\sf}

\title{The second minimum of Barnes-Wall lattices}
\author{Gabriele Nebe\footnote{nebe@math.rwth-aachen.de}}  
\date{Lehrstuhl f\"ur Algebra und Zahlentheorie, 
RWTH Aachen University, Germany}
 
\begin{document}
\maketitle

 \begin{center} 
	 {\em In memory of Jacques Martinet} 
 \end{center} 

 {\sc Abstract.} 
 The paper gives a recursive construction of the Barnes-Wall lattices 
 as subdirect products. This is used to show that 
 the Barnes-Wall lattices of minimum $d$ do not contain any vectors of norm 
 $a$ with $d<a<3d/2$.
 \\
 On donne une construction recursive des r\'eseaux de Barnes-Wall qui permet de 
 d\'emontrer que 
 les r\'eseaux de Barnes-Wall de minimum $d$ ne poss\`edent aucun vecteur de norme 
 $a$ avec $d<a<3d/2$.
  \\
  {\sc keywords:}  Barnes-Wall lattices, minima of lattices
  \\

\section{Introduction} 

The Barnes-Wall lattices form an infinite series of lattices 
of dimension $2^m$ for $m\in \N $. 
They are constructed in 
\cite{BW} which explicitly elaborates a subset of their minimal vectors 
 to show that the Barnes-Wall lattices are locally densest lattices.  
This fact also follows from an inspection of their automorphism group.
 For $m\neq 3$ this group is a subgroup of index 2 of the real Clifford group 
 (see \cite{NRS}), a fact that allows  Bachoc in \cite{Bachoc} 
 to show that for $m\geq 3$ all non-empty
 layers of the Barnes-Wall lattices form spherical 6-designs. 
 This paper gives a construction of the Barnes-Wall lattices as 
 subdirect products (Theorem \ref{subdir}) that is used to give an easy proof 
 of their kissing number (Lemma \ref{lemmin}) 
 and to show that for $m\geq 3$ the minimal vectors
 form a spherical 6-design (Proposition \ref{des}). 

 The most prominent construction of the Barnes-Wall lattices is by applying 
 Construction D to a chain of Reed-Muller codes (see for instance 
 \cite[Chapter 8, Section 8]{SPLAG}, \cite{BW}, 
 and also \cite{HN} for a basis-independent formulation). 
 Berlekamp and Sloane \cite{RM} show that 
 in the $r$-th order binary Reed-Muller code of length $2^m$
 and minimum distance $d=2^{m+r}$, the only codewords having weight between $d$ and $2d$
 are those with weights of the form $2d-2^i$ for some $i$. 
 Motivated by this observation, certain experiments, and the theta series 
 of the Barnes-Wall lattice of dimension 64 and 128 in \cite{OEIS}, 
 Christoph Keller \cite{Keller} conjectured  that a similar property 
 should also be true for the Barnes-Wall lattices.
 
 Theorem \ref{secondmin} is a first step in this direction showing that 
 the Barnes-Wall lattices of minimum $d$ have no vectors of norm $a$ with 
 $d<a <3d/2 $.

\section{The two Barnes-Wall lattices of dimension $2^m$} 

In \cite{BW} Barnes and Wall construct a series of lattices in dimension 
$N:=2^m$ for any $m\in \N$. 
Put 
 $\cV_m = \F_2^m$ to denote the $m$-dimensional vector space over the field with 
2 elements and fix a basis $(v_1,\ldots , v_m)$ of $\cV_m $. 
Put 
 $${\mathcal T}_r(m):= \{ \cU \leq \cV_m \mid \cU = \langle v_i \mid i\in I \rangle _{\F _2},
I \subseteq \{ 1,\ldots , m \},  |I| = r \} $$
to denote the set of $r$-dimensional subspaces of $\cV_m$ that have a basis 
that is a subset of $\{ v_1,\ldots , v_m \}$. 

We now let $\cV_m $ index the elements of an orthonormal basis 
 $$(e_v \mid v\in \cV_m )$$ of the euclidean space $(\R^N, (\phantom{x},\phantom{x}))$. 
 Let 
$$\Gamma_m := \langle e_v \mid v\in \cV_m \rangle _{\Z}  $$
denote the standard lattice spanned by this orthonormal basis. 
For any subset $\cU $ of $\cV_m $ we put
$$x_{\cU} := \sum _{v\in \cU} e_v \in \Gamma_m .$$
Then \cite{BW} defines the following sublattices of $\Gamma_m $:

\begin{definition}\label{defbw}
Let $\lambda := (\lambda_0,\ldots , \lambda _m) \in \Z^{m+1} $
        such that 
 $\lambda _0 = 0, \lambda _r -1 \leq \lambda _{r-1} \leq \lambda _r$
for all $1\leq r\leq m$.
Then 
$$\Lambda (\lambda ):= \langle 2^{\lambda _{m-r}} x_{\cU} \mid 
 0\leq r \leq m, \cU \mbox{ affine subspace of } \cV_m , \dim(\cU) =  r \rangle _{\Z} .$$
\end{definition}

For $\lambda $ as in Definition \ref{defbw} 
        put $\lambda ' := (\lambda _0',\ldots , \lambda _m') $ where
        $\lambda _r' = \lambda _m - \lambda _{m-r} $. 
Then \cite[Theorem 3.1, Theorem 3.2]{BW} give the following properties of the lattices $\Lambda (\lambda) $.

\begin{proposition}\label{propBW}
        \begin{itemize}
                \item[(a)]
                $\Lambda (\lambda ') = 2^{\lambda _m} \Lambda (\lambda )^{\# }$.
                \item[(b)]
                        A $\Z $-basis of $\Lambda (\lambda )$ is given by
			$$\bigcup_{r=0}^m \{ 2^{\lambda _{m-r}} x_{\cU } \mid \cU \in {\mathcal T}_r (m) \} .$$
               \item[(c)]
$\det(\Lambda (\lambda ) ) = 
                        2^{2d}$ where  $d=\sum_{r=0}^m \lambda _r {{m}\choose{r}} $.
\item[(d)]  $\min(\Lambda (\lambda )) = 2^{\alpha }$ where
$\alpha = \min \{ m-r+2\lambda _r \mid 0\leq r \leq m \} $.
        \end{itemize}
\end{proposition}

Barnes and Wall single out two particular lattices among the
lattices $\Lambda (\lambda ) $:

\begin{theorem} \label{BW} (\cite[Theorem 4.3]{BW}) 
Put $\lambda _r := \lfloor \frac{r}{2} \rfloor $ and 
 $\mu _r := \lfloor \frac{r+1}{2} \rfloor $ ($0\leq r \leq m $) and put 
	$$\Lambda _m := \Lambda (\lambda ), \ \Delta _m := \Lambda (\mu ) .$$ 
	The index of $\Lambda _m$ in the standard lattice 
	$\Gamma _m $ is $2^j$ where
	$$j= \log_2([\Gamma _m : \Lambda _m ]) = \sum _{r=0}^{m} \lfloor \frac{r}{2} \rfloor 
	{{m}\choose{r}} = (m-1) 2^{m-2} .$$
Then $\Lambda_m \supset \Delta _m \supset 2 \Lambda _m $ 
with $[\Lambda _m : \Delta _m] = [\Delta _m : 2\Lambda _m ] =2^{m-1}$. 
	Moreover $\Lambda_m $ and $\Delta _m$ are similar lattices. 
For the minimum of the two lattices 
	we get $\min (\Lambda_m ) = 2^{m-1}, \min (\Delta_m ) = 2^m$. 
\end{theorem} 

The lattices that are commonly known as ``the'' Barnes-Wall lattices are 
scaled versions of the lattices $\Lambda _m$ from Theorem \ref{BW}.

\begin{remark}\label{rescaleBW}
If $m\geq 3$ is odd then $\BW_m :=2^{-(m-1)/4} \Lambda _m$ is an even unimodular lattice 
	of minimum $\min (\BW_m) = 2^{(m-1)/2} $. 
	\\
If $m$ is even then  $\BW_m:=2^{-(m-2)/4} \Lambda _m$ is an even $2$-modular lattice 
	of minimum $\min(\BW_m) = 2^{m/2} $. 
\end{remark}

The automorphism group $G_m\cong 2^{1+2m}_+.\Omega _{2m}(2)$ of 
$\Lambda _m$ (for $m=3$ we put $G_m$ to denote the
stabiliser of $\Delta_m $ in $\Aut(\Lambda _m)$) 
is a normal subgroup of index 2 in the real Clifford group 
$\langle G_m , h \rangle \leq \GL_N(\R)$ (see for instance \cite{Bachoc}, \cite{NRS}).  
The element $\sqrt{2} h$ is rational and induces the similarity 
between 
$\Lambda _m $ and $\Delta _m= \sqrt{2}h \Lambda _m$.
The group $G_m$ is studied in detail in \cite{BE}. 

\begin{theorem} (\cite[Th\'eor\`eme II.4]{BE}) \label{trans} 
$G_m$ acts transitively on the set of minimal vectors of $\Lambda _m$.
\end{theorem}

\section{The second minimum of the lattices $\Lambda_m $.} 

\subsection{The minimal classes in $\Lambda_m/\Delta _m $}

The key observation for having a recursive proof of the second minimum 
of the Barnes-Wall lattices is given in the following lemma. 

\begin{lemma} \label{minclass} 
Let $\ell \in \Lambda_m$ be a minimal vector, i.e. $(\ell ,\ell ) = 2^{m-1}$. 
Then for any $x\in \ell +\Delta _m $ either $(x,x) = (\ell ,\ell ) = 2^{m-1}$ or 
$(x,x) \geq 2^m$. 
\end{lemma} 

\bew
By Theorem \ref{trans} the group $G_m$ acts transitively on the set of minimal 
classes of $\Lambda_m/\Delta _m$ so we make a suitable choice of the
minimal vector $\ell $. 
\\
First assume that $m$ is odd. Then we choose the minimal vector 
$\ell = 2^{(m-1)/2} e_0 \in \Lambda _m$ and let $x = \ell +d \in \ell +\Delta _m$ 
with $(x,x) > 2^{m-1}$.  
Write $x = \sum _{v\in \cV_m } a_v e_v $ with coefficients $a_v\in \Z $ 
in the orthonormal basis and let $2^j$ be the maximal $2$-power that 
divides all $a_v$. If $j \geq (m-1)/2$ then 
$(x,x)$ is a multiple of $2^{m-1}$ and hence $(x,x) \geq 2^m$ by the assumption that 
$(x,x) > 2^{m-1}$. 
\\
So assume that $j < (m-1)/2$ and put 
$$y:= 2^{-j} x = \sum _{v\in \cV_m } b_v e_v $$ with $b_v = 2^{-j} a_v \in \Z $. 
Then  $2^{-j} d = (b_0-2^{(m-1)/2-j})e_0 + \sum_{0\neq v\in \cV_m } b_v e_v $
and the set 
$S:=\{ v\in \cV_m \mid b_v \mbox{ odd } \} $ is the set of indices of the odd coefficients 
in $2^{-j} d$. By \cite[Lemma 3.3]{BW} the cardinality of $S$ is $\geq 2^{m-2j}$ and hence 
$(x,x) = 2^{2j} (y,y) \geq 2^{2j} |S| \geq 2^m $. 
This shows the lemma if $m$ is odd. 
\\
For even $m$, there are no minimal vectors of $\Lambda _m$
that are scalar multiples of one of the $e_v$. 
However, as $\Lambda _m$ and $\Delta _m$ are similar, 
we may use the same argument as before for  the minimal classes of 
$\Delta _m/ 2\Lambda _m$. 
So we choose the minimal vector 
 $\ell  = 2^{m/2} e_0 \in \Delta _m $ and assume that $d\in 2\Lambda _m$ 
is such that 
$x=\ell +d\in \ell +2\Lambda _m$ has norm $(x,x) > 2^m = \min (\Delta _m)$. 
Then the same argument as before shows that $(x,x) \geq 2^{m+1}$ which 
shows the lemma also for $m$ even. 
\eb

\subsection{A recursive construction of $\Lambda _m$ as a subdirect product}

The Barnes-Wall lattices have an easy construction as a subdirect product, 
very similar to the doubling construction for the Reed-Muller codes. 
Fixing the basis $(v_1,\ldots , v_{m+1}) $ of $\cV _{m+1}$ as before, 
there is a natural embedding 
$$\iota : \cV_m  = \langle v_1,\ldots , v_m \rangle \hookrightarrow \cV_{m+1} , \iota (v) = v  .$$ 
Combining $\iota $ with the translation along $v_{m+1}$ we obtain a bijection 
$$\tau : \cV_m \to v_{m+1} + \iota (\cV_m) , v \mapsto v_{m+1} + \iota(v) $$ 
so that $\cV_{m+1} = \iota (\cV_m) \disj \tau (\cV_m )$. 
By abuse of notation we also use $\iota $ and $\tau $ to denote the $\Z $-linear maps 
$$ \iota : \Gamma _m \to \Gamma _{m+1} , e_v \mapsto e_{\iota(v)} 
\mbox{ for all } v\in \cV_m $$ 
and 
$$\tau : \Gamma _m \to \Gamma _{m+1} , e_v \mapsto e_{\tau (v)} 
\mbox{ for all } v\in \cV_m .$$ 
Then $\iota $ and $\tau $ are isometric embeddings
and $\Gamma _{m+1}$ is the orthogonal sum of $\iota (\Gamma _m) $ and $\tau (\Gamma _m)$.
In this notation we obtain

\begin{theorem}\label{subdir}
	$\Lambda _{m+1} = \{ \iota(\ell ) + \iota (d) + \tau (\ell ) 
	\mid \ell \in \Lambda _m, d\in \Delta _m \} $. 
\end{theorem} 

\bew
The right hand side, $\Lambda $, is a lattice, in fact we have
$$\Lambda = (\iota + \tau ) (\Lambda _m) \oplus \iota (\Delta _m)  .$$
So the index of $\Lambda $ in $\Gamma _{m+1}$ is $2^j$ with
$$
\begin{array}{l} j = \log_2([\Gamma _{m+1} : \Lambda ]) = 
\log_2( [\Gamma _m : \Delta _m]) + \log_2( [\Gamma _m : \Lambda _m ] )  = \\
(m+1) 2^{m-2} + (m-1) 2^{m-2}   = m 2^{m-1} = \log_2 ([\Gamma _{m+1} : \Lambda _{m+1}] ). \end{array} $$
So it remains to show that the 
basis of $\Lambda _{m+1} $ given in Proposition \ref{propBW} (b) is contained in $\Lambda $. 
So let $\cU \in {\mathcal T}_r(m+1)$. 
\\
If $v_{m+1} \in \cU $ then $\cU' := \cU \cap \cV_m \in {\mathcal T}_{r-1}(m) $
and 
$$2^{\lfloor \frac{m+1-r}{2} \rfloor }x_{\cU} = \iota (2^{\lfloor \frac{m-(r-1)}{2} \rfloor }x_{\cU'} ) + \tau (2^{\lfloor \frac{m-(r-1)}{2} \rfloor } x_{\cU '}) 
\in (\iota +\tau) (\Lambda _m) .$$
If $v_{m+1} \not\in \cU $ then $\cU \in {\mathcal T}_{r}(m) $
and 
$$2^{\lfloor \frac{m+1-r}{2} \rfloor }x_{\cU} = \iota (2^{\lfloor \frac{m-r+1)}{2} \rfloor }x_{\cU'} ) \in \iota (\Delta _m ).$$
\eb

For related constructions see for instance \cite{Griess} and 
\cite{NRSBW}.

\subsection{The main result} 

\begin{theorem} \label{secondmin}
Let $x\in \Lambda_m $ be such that $(x,x) > \min (\Lambda _m) = 2^{m-1}$. Then $(x,x) \geq 2^{m-1} + 2^{m-2}$. 
\end{theorem} 

\bew
We proceed by induction on $m$. 
The cases $m=2,\ldots,5$ follow immediately from Remark \ref{rescaleBW}, 
as here $\BW_m$ is an even lattice of minimum $\leq 4$.

For the induction step assume that the theorem holds for $m$ and let 
$x = \iota(\ell ) + \iota (d) + \tau (\ell )  \in \Lambda _{m+1}$ with 
$$2^m < (x,x) = (\ell +d , \ell +d ) + (\ell , \ell) < 2^{m} + 2^{m-1} .$$
Without loss of generality we assume that $(\ell, \ell) \leq (\ell+d,\ell+d)$. 
\\
Then there are three cases to consider: \\
If $\ell = 0$, then $d\in \Delta _m$ is not a minimal vector as $(x,x) = (d,d) > 2^m$. 
Hence by induction hypothesis (using the fact that $\Delta _m $ is similar to $\Lambda _m$) we have $(x,x) = (d,d) \geq 2^{m} + 2^{m-1}$. 
\\
If $\ell\in \Lambda _m$ is not a minimal vector 
then by assumption $(\ell,\ell) \geq 2^{m-1}+2^{m-2}$ and hence 
$(x,x) \geq 2(\ell,\ell) \geq 2^m+2^{m-1}$. 
\\
In the last case 
 $\ell \in \Lambda _m$ is a minimal vector, i.e. $(\ell,\ell) = 2^{m-1}$.
By assumption $(x,x) > 2^m$, therefore $\ell+d\in \ell+\Delta _m $ is 
not a minimal vector, so $(\ell+d,\ell+d) \geq 2^{m} $ by Lemma \ref{minclass}. 
Therefore $(x,x) \geq 2^{m} + 2^{m-1}$. 
\\
Combining these three cases shows  the claim for 
$m+1$ and finishes the induction step. 
\eb

\subsection{Some further consequences} 

Theorem \ref{subdir} also allows to deduce an easy recursive formula for the set
$$S(\Lambda _m):= \{ \ell \in \Lambda _m \mid (\ell , \ell) = 2^{m-1} \} $$ 
of minimal vectors in $\Lambda _m$ and its cardinality $s_m := | S(\Lambda _m) |=|S(\Delta _m) |$.

\begin{lemma} \label{lemmin}
	$$S(\Lambda _{m+1}) = \begin{array}{l} 
\{ (\iota + \tau )(\ell)  + \iota (d) \mid \ell , \ell + d \in S(\Lambda_{m}), 
d\in \Delta _{m} \} \cup \\
\{ \iota(d) \mid d\in S(\Delta _{m}) \} \cup 
\{ \tau (d) \mid d\in S(\Delta _{m}) \} \end{array} .$$
	The kissing number $s_{m+1}$ of $\Lambda _{m+1}$ satisfies 
	$s_1 = 4$ and for $m\geq 1$ 
	$$s_{m+1} = (2^{m+1}+2) s_{m} .$$
\end{lemma}

\bew
For any two minimal vectors  $\ell_1,\ell_2 \in \Lambda _m $ with $(\ell_i,\ell_i ) = 2^{m-1}$ and $\ell_1 + \Delta _m = \ell _2 + \Delta _m$ we have that 
$\ell _1 \pm \ell _2 \in \Delta _m$ and hence 
$2^{m} \pm 2 (\ell_1,\ell_2) \geq 2^m 
= \min (\Delta _m) .$
So $(\ell _1,\ell_2) =0$ and the class $\ell _1 + \Delta _m$ contains at most 
$2^{m+1}$  minimal vectors, 
$ \pm \ell_1, \pm \ell _2, \ldots , \pm \ell _{2^m}   $ with
 $(\ell _i, \ell _j) = 0$ for $i\neq j$. 
 To show that there are at least $2^{m+1}$ such minimal vectors in 
 each class in $\Lambda _m / \Delta _m$ and hence also in $\Delta _m / 2\Lambda _m$
 we use Theorem \ref{trans} and the fact that $\Lambda _m$ and $\Delta _m$ 
 are similar lattices to choose a suitable class. 
 If $m$ is odd put  $\ell +\Delta _m=  2^{(m-1)/2} e_0  +\Delta _m$ 
 and for even $m$ choose 
 the class $\ell + 2\Lambda _m = 2^{m/2} e_0  +2\Lambda _m \in \Delta _m/2\Lambda _m$.
 These classes contain the $2^{m+1}$ 
 minimal vectors 
 $\pm 2^{\lfloor m/2 \rfloor } e_v$ for all $v\in \F_2^m$.
 \\
 So there are exactly
 $2^{m+1} s_{m}$ minimal vectors of the form 
$(\iota + \tau )(\ell)  + \iota (d) \in \Lambda _{m+1}$, where $\ell , \ell+d  \in \Lambda _{m}$ 
are both minimal vectors in the class $\ell + \Delta _{m}$ and 
$2s_{m}$ minimal vectors of the form $\iota (d)$ or $\tau (d)$ where 
$d\in \Delta _{m}$ is a minimal vector. 
\eb

Using this description of the minimal vectors of $\Lambda _m$ 
it is easy to show that for $m\geq 3$ 
these minimal vectors form a spherical $6$-design. 
Though this is a well known fact (see for instance \cite[Corollary 5.3]{Bachoc}) 
the direct proof given below is substantially easier than the one using the description in
\cite[Th\'eor\`eme II.2]{BE} of the minimal vectors of the Barnes-Wall lattices. 

Let $X \subset \R ^N$ be a finite subset of the sphere of squared radius 
$M$ such that 
such that $-x \in X $  for all $x\in X$.
By \cite[Th\'eor\`eme 3.2, (3.6)]{Venkov} the set 
$X$ is a spherical $6$-design, if and only if 
$$ \sum_{x\in X}  (x,v)^6 = |X| M^3 (v,v)^3 \frac{1\cdot 3 \cdot 5}{N(N+2)(N+4)} \mbox{ for all } v\in \R^N .$$
For $z\in \R $ we put
$$c_m(z) := 
 s_m \min(\Lambda_m)^3 z^3 \frac{1\cdot 3 \cdot 5}{2^m(2^m+2)(2^m+4)}  .$$
 Then $S(\Lambda_m)$ is a spherical $6$-design, if and only if 
 for all $v\in \R^{2^m}$ the sum
$\sum _{x\in S(\Lambda _m)} (x,v)^6 = c_m((v,v)) $.
As $\min (\Delta _m) = 2 \min (\Lambda _m)$ the set 
 $S(\Delta_m)$ is a spherical $6$-design, if and only if 
$\sum _{x\in S(\Delta _m)} (x,v)^6 = 8 c_m((v,v)) $.
Using induction over $m$ the construction in Theorem \ref{subdir} allows us to
give a quick proof of the following well known result.

\begin{proposition}  \label{des}
For $m\geq 3$ the sets 
$S(\Lambda _m) $ and $S(\Delta _m)$ form spherical $6$-designs. 
\end{proposition} 

\bew
We use the criterion in \cite[Th\'eor\`eme 8.1]{Venkov} and show that 
$$\sum _{x,y \in S(\Lambda _m) } (x,y) ^6 = s_m c_m(\min (\Lambda _m)) .$$
Then the fact that also $S(\Delta _m)$ is a spherical $6$-design follows from 
the similarity of these two lattices. 
\\
By Theorem \ref{trans} 
the automorphism group of $\Lambda _m$ acts transitively on $S(\Lambda _m)$ 
so we may choose a suitable $x_0 \in S(\Lambda _m)$ and 
prove that 
$$\sum_{x\in S(\Lambda _m)}  (x,x_0)^6 = c_m(\min (\Lambda _m)) .$$
We proceed by induction on $m$, the case $m=3$ can be settled by direct computations. 
Assume that the proposition is proven for $m$ and choose 
$x_0 = \tau (d_0) \in S(\Lambda _{m+1}) $ for some $d_0\in S(\Delta _m)$. 
Then by the description of $S(\Lambda _{m+1})$ in Lemma \ref{lemmin} 
$$(\star) \ \ \sum_{x\in S(\Lambda _{m+1})}  (x,x_0)^6 = 2^{m+1}  \sum _{\ell \in S(\Lambda_m)} 
(\ell , d_0)^6 + \sum _{d\in S(\Delta _m)} (d,d_0)^6 .$$
By induction hypotheses $S(\Lambda _m)$ and $S(\Delta _m)$ both form 
spherical $6$-designs and hence 
$$(\star)  = (2^{m+1}+8) c_m(\min(\Delta _m)) =  
2^3(2^{m-2}+1)  s_m 2^{6m-3}  \frac{15}{2^m(2^m+2)(2^m+4)}  =  $$
$$ s_{m+1} 2^{6m} \frac{15(2^{m-2}+1)}{2^{m+3}(2^{m-1}+1)(2^{m-2}+1)(2^m+1)} =
c_{m+1}(\min(\Lambda_{m+1})).$$
\eb

\end{document}